\documentclass[noinfoline, 11pt]{imsart}
\usepackage{amssymb,amsfonts,amsmath,latexsym,mathtools,enumerate}
\usepackage[table]{xcolor}
\usepackage{url}
\usepackage{hyperref}
\definecolor{burgundy}{rgb}{0.5, 0.0, 0.13}
\definecolor{camel}{rgb}{0.76, 0.6, 0.42}
\definecolor{chamoisee}{rgb}{0.63, 0.47, 0.35}
\definecolor{grey1}{RGB}{128,128,128}

\hypersetup{colorlinks=true,linkcolor=burgundy,citecolor=teal,urlcolor=burgundy,linktoc=page}

\newtheorem{theorem}{Theorem}
\newtheorem{lemma}[theorem]{Lemma}

\newcommand\Z{{\mathcal Z}}
\newcommand\N{{\mathbb N}}
\newcommand\R{{\mathbb R}}

\newcommand\cL{{\mathcal L}}
\newcommand\balpha{\boldsymbol\alpha}

\newcommand\E{{\mathbb E}}

\newcommand{\Cov}{\mbox{\rm Cov}}
\newcommand{\Var}{\mbox{\rm Var}}
\newcommand{\vol}{\mbox{\rm vol}}
\newcommand{\sinc}{\mbox{\rm sinc}}
\newcommand{\cs}{\mbox{\rm cs}\ }

\begin{document}
\begin{frontmatter}
\title{A note on $3d$-monochromatic random waves and cancellation}
\date{\today}
\begin{aug}
\author{\fnms{Federico} \snm{Dalmao}\ead[label=e1]{fdalmao@unorte.edu.uy}}
\thanks{
DMEL, CENUR Litoral Norte, 
Universidad de la Rep\'{u}blica, Uruguay. 
fdalmao@litoralnorte.udelar.edu.uy}
\end{aug}

\begin{abstract}
In this note we prove that the asymptotic variance of the nodal length of complex-valued monochromatic random waves restricted to an increasing domain in $\R^3$ is linear in the volume of the domain. 
Put together with previous results this shows that a Central Limit Theorem holds true for $3$-dimensional monochromatic random waves. 
We compare with the variance of the nodal length of the real-valued $2$-dimensional monochromatic random waves where a faster divergence rate is observed, 
this fact is connected with Berry's cancellation phenomenon. 
Moreover, we show that a concentration phenomenon takes place.

\end{abstract}
\end{frontmatter}

\thispagestyle{empty}
\maketitle

\section{Introduction}
During the last decade an impressive effort has been dedicated to understanding the behavior of geometric functionals of 
random waves in the high energy limit. 
Motivation often comes from Quantum Chaos as it is concerned with the behavior of Laplace eigenfunctions 
as the frequency or the energy level (equivalently,  the eigenvalue) tends to 
infinity. For a general picture of the field, see the seminal papers \cite{Be02,BD00}, the recent reviews \cite{wigmanr,marinuccir,vidottor} 
and the background section in \cite[S.1.6]{kkw} for an example of an early stage discussion on this topic. 

Some instances of geometric functionals are the length of nodal (i.e: zero) curves, the number of 
critical points and the number of connected components of the zero sets, see e.g. \cite{BD00}. 
Within this general framework, we focus our attention on the nodal length of complex-valued random waves defined on $\R^3$ and 
compare it to the well studied case of the nodal length of real-valued random waves defined on $\R^2$. 
The analogy between these two cases is apparent when looking at the covariance functions, the analytic expressions and the expansions of the 
nodal lengths.
Indeed, since we use real analytic methods we can think the former random waves as $\R^2$-valued. 
See sections \ref{s:prelim}, \ref{s:proof} and \ref{s:anciliary}. 
\smallskip

A central family of random fields (either real or complex-valued, defined on $\R^d$, $d\geq 3$), 
representing monochromatic waves, 
is the so-called Berry's Random Wave Model (RWM for short). 
M. Berry \cite{Be77} conjectured that the RWM can serve as a model 
for deterministic '{}generic'{} Laplace eigenfunctions on manifolds with negative 
curvature in the high energy regime. 
Furthermore, 
the RWM can be thought of as a universal model since 
it can be proved that, in some local sense, 
it is the limit of other important Laplace random eigenfunctions models 
as Random Spherical Harmonics (RSH) \cite[Eq.(1.8)-(1.9)]{wigmanf}
or Arithmetic Random Waves (ARW) defined in the flat torus \cite[S1.6.1]{kkw}. 
See also \cite{CH16} and \cite{gass} for the case of Riemannian manifolds and \cite{rs} for some deterministic related constructions. 
For a discussion on these approximations, see section 1.6. in \cite{kkw}. 
\smallskip

The complex (resp. real)-valued RWM can be defined as the complex (resp. real)-valued centered stationary isotropic Gaussian random field $\psi$ defined on $\R^d$ having covariance function
\begin{equation} \label{d:r-gral}
 r(x,y) = \E[\psi(x)\overline{\psi(y)}] 
 = c_d \frac{J_{\lambda}(|x-y|)}{|x-y|^\lambda},
\end{equation}
where $c_d$ is a normalizing constant, $|\cdot|$ stands for the Euclidean norm in $\R^d$ and $J_\lambda$ is the Bessel function of the first kind with index $\lambda=\frac{d-2}{2}$.
Alternatively, the RWM can be defined by stating that its spectral measure is uniform on the unit sphere. 
In the real case, the conjugation in \eqref{d:r-gral} shall be omitted.
\smallskip

More generally, set $\psi_k(x) = \psi(kx):k\ge 1$. 
Note that $\psi=\psi_1$. 
Then, it is well known that for $k\geq 1$, $\psi_k$ is a solution of the Helmholtz equation
\begin{equation*}
 \Delta \psi_k (x) + k^2 \psi_k (x) = 0, \ x\in\R^d, 
\end{equation*}
which explains the terminologies frequency, energy or eigenvalue used to refer to $k^2$. 
See section 1.1. in \cite{npr}.

It is worth mentioning that the study of the nodal sets of $\psi$ on growing domains is equivalent to the study of the nodal sets of $\psi_k$ 
restricted to a fixed domain  
in the high energy limit, i.e: as $k\to\infty$, 
see remark 3.6. in \cite{del} and remark 1.3 in \cite{npr}. 
In this note we choose to work on the former framework 
because our natural domain is the non-compact manifold $\R^3$. 
See the discussion in section \ref{s:discussion} below.
\smallskip

The nodal sets of these models were deeply studied in the last years, though the major part of the literature focuses on the 
real-valued $2$-dimensional 
case. The papers \cite{npr} (concerning real and complex-valued $2$ dimensional RWM) and 
\cite{del} (concerning complex-valued $3$-dimensional RWM) are key for our work. 
See also the recent works \cite{pv}, \cite{npv} dealing with functional convergence on the $2$-dimensional case
and the reviews \cite{marinuccir,vidottor,wigmanr}. 
\smallskip

In the real-valued $2$-dimensional case the RWM 
soon showed a surprising behavior concerning the high energy limit of the variance of their nodal 
length (i.e: the length of their zero curves) : it diverges far more slowly than anticipated 
(similar results hold true for some related models as RSH and ARW). 
More precisely, while the expectation and the variance of the 
nodal length were expected to be of the same order as $k\to\infty$, 
rephrasing equations (1.7)-(1.8) in \cite{npr}, 
it happens that
\begin{equation*}
 \Var(\cL_k(B)) \sim\ \log(\E[\cL_k(B)]),\ {\rm as}\ k\to\infty, 
\end{equation*}
where $B$ is a fixed domain in $\R^2$ and $\cL_k(B)$ is the nodal length associated 
to $\psi_k$ restricted to $B$. 
This fact was first predicted by Berry \cite{Be02} 
where he expressed that this behavior is due to 'an obscure cancellation phenomenon'. 
\smallskip

The first rigorous confirmations of the so called 
Berry's cancellation arrived in \cite{wigmanf} in the case of RSH, 
in \cite{mprw} for ARW in the flat torus 
and in \cite{npr} for the planar RWM. 
One important feature of the last reference is that it was the first one to deal with non-compact manifolds, 
which is also the case we are interested in.

The early results, obtained using Rice formulas, did not allow a clear understanding of the cancellation phenomenon, 
see e.g. \cite{Be02,wigmanf,kkw}. 
The use of Wiener Chaos techniques, see e.g. \cite{mprw,npr,DNPR17}, provided a precise explanation 
of the cancellation phenomenon in terms of the Wiener (a.k.a. chaotic or Hermite) expansion of the nodal 
length which combines in a subtle manner geometric and probabilistic aspects. 
Roughly speaking, the nodal length restricted to a Borel set $B$, $\cL(B)$ say, 
can be expanded in the $L^2$-sense in the form
\begin{equation*}
\cL(B)  = \sum^\infty_{q=0} I_{2q}(B),
\end{equation*}
where the components $I_{2q}(B)$ are orthogonal random variables. 
In particular, $I_0(B)$ equals the expectation of $\cL(B)$ and 
the rest of the terms are combinations of integrals involving polynomials 
on the underlying process and its derivative. 
It was observed that the second component $I_2(B)$ vanishes due to geometric reasons 
reflected in the coefficients of the expansion, 
see \cite{notarnicola} for a general treatment of this issue. 
The cancellation of the coefficients in the second term in the expansion avoids the integrals 
appearing in the computation of the variance of $I_2(B)$
to produce higher rates of divergence, 
as those which do take place in the case of non-zero levels, 
see e.g. \cite{mprw}. 
See section \ref{s:discussion} below.
\medskip

In \cite{del} complex-valued $3$-dimensional Berry's general model in growing domains was considered. 
In particular, for monochromatic random waves 
(verifying \eqref{d:r-gral} with $d=3$ 
or, equivalently, verifying \eqref{d:r} below)  
it was shown that the variance of the second component in the Wiener expansion vanishes asymptotically and that the variance of $\cL(B)$ grows at most linearly in the volume $\vol(B)$ as the underlying domain $B\uparrow\R^3$. 
Nevertheless, no lower bound was pursued. 
Consequently, the possibility of having a strictly lower order variance remained open. 
In view of the behavior of the variance in the real-valued $2$-dimensional case, the question about the 
true order of the variance 
arises. 
\medskip

In the present note, 
we establish that the variance of the nodal length of the complex-valued  $3$-dimensional RWM 
is linear (with a strictly positive coefficient) in the volume of the growing domain, 
which coincides with the order of the expectation. 
This behavior differs from the observed one in the real valued $2$-dimensional case \cite{npr} where 
the variance diverges faster than the expectation, 
\begin{equation*}
 \E[\cL(B)] \approx {\rm area}(B);\quad  
 \Var(\cL(B)) \approx\ {\rm area}(B) \log({\rm area}(B)),\ {\rm as}\ B\uparrow\R^2. 
\end{equation*}
A rough explanation may be as follows. 
In both cases the second component $\Var(I_2(B))$ vanishes due to geometric reasons as said above. 
The key difference between these cases lies in the integrability of the covariance functions 
which is a fact of probabilistic nature (short vs. long memory, say). 
In the $2$-dimensional case $r(x,y) = c_2 J_0(|x-y|)$ belongs to $L^6(\R^2)$ but not to $L^4(\R^2)$ 
while in the $3$-dimensional case $r(x,y) = J_{1/2}(|x-y|)/|x-y|^{1/2}$ belongs to $L^4(\R^3)$. 
Grosso modo, see lemma 9 in \cite{del}, 
once normalized by the area or volume, according to the dimension, 
the variance of $I_{2q}(B)$ is written as an integral of the $2q$-th power of $r$ 
(actually, the product of $2q$ factors chosen from $r$ and its derivatives). 
Thus, the normalized variance of $I_4(B)$ diverges in dimension $2$ while the rest of the variances (i.e: 
those of $I_{2q}(B):q\geq 3$ in $2$d, those of $I_{2q}(B):q\geq 2$ in $3$d), once normalized, converge. 
Observe that it follows that the variance is of lower order than the square of the expectation as $B$ increases, 
implying the concentration of the probability. 
See section \ref{s:discussion} for some details.
\medskip
 
To prove our result, we carry on a careful analysis of the variance of the fourth component of the Wiener expansion 
of $\cL(B)$ which involves a few dozen terms. 
We need to compute explicitly each coefficient in the expansion and to find the asymptotic of the involved integrals. 
Since we look for a lower bound for the variance, we can omit a convenient number of nonnegative terms.
\medskip

Let us fix some notation: $\cs$ denotes an unimportant constant whose value may change from one line to another, 
$a_R\approx b_R$ means that $\lim_{R\to\infty}\frac{a_R}{b_R}=c$ with $c>0$ and  
$a_R\sim b_R$ when $\lim_{R\to\infty}\frac{a_R}{B_R}=1$, we write $a_R\gtrsim b_R$ when there exists $B'_R\sim b_R$ 
such that $a_R\geq B'_R$.

\section{Aknowledgements}
Je voudrais remercier Giovanni Peccati et Jos\'e R. Le\'on pour les discussions sur ce sujet.

\noindent The author was partially supported by Agencia Nacional de Investigaci\'on e Innovaci\'on (ANII). 

\section{Problem setting and main result}
Define the $3$-dimensional monochromatic random wave model (RWM) as the complex-valued centered stationary isotropic Gaussian random field 
\begin{equation} \label{d:psi}
\psi(x) := \xi(x) + i\eta(x),\quad x\in\R^3, 
\end{equation}
such that its real and imaginary parts $\xi,\eta$ are centered independent Gaussian random fields with 
covariance function
\begin{equation} \label{d:r}
r(x) := \sinc(|x|),
\end{equation}
where $\sinc$ stands for the cardinal sine function $\sinc(\cdot)=\frac{\sin(\cdot)}{\cdot}$ and $|\cdot|$ is the Euclidean norm in $\R^3$. 
Observe that this coincides with \eqref{d:r} for $d=3$.
\smallskip

For any bounded domain $B$ in $\R^3$, we introduce the nodal curve $\Z(B)$ 
and the nodal length $\cL (B)$ by
\begin{align*} 
\Z(B) &:=\big\{x\in B\,:\,|\psi(x)|=0\big\}, \notag\\
\cL(B) &:={\rm length}(\Z(B)).
\end{align*}
In \cite{del} it was proved that $\E[\cL(B)]=\frac{1}{3\pi}\vol(B)$ (Corollary 2 with $\lambda=\frac13$), 
that $\Var(\cL(B)) = O(\vol(B))$ as $B\uparrow\R^3$ and that a Central Limit Theorem takes place (proposition 10). 
Since the limit variance was not known to be positive, the limit law could be degenerated.
\medskip

The main result of the present note is contained in the next theorem. 
\begin{theorem} \label{t:main}
 Let $\psi$ be defined as in \eqref{d:psi}-\eqref{d:r}.
 Consider $B_R$ the centered ball of radius $R$ in $\R^3$, then  
 there exists a constant $c>0$ such that 
 \begin{equation*}
 \lim_{R\to\infty} \frac{\Var(\cL(B_R))}{\vol(B_R)}\geq c.
 \end{equation*} 
\end{theorem}
\bigskip

Observe that in \cite{del} the underlying domain was $Q_n=[-n,n]^3$ as $n\to\infty$ while here we use $B_R$ as $R\to\infty$. We choose the latter for simplicity of exposition, but one can easily adapt the computation of the expectation and the variance from $B_R$ to $Q_n$ and conversely.  
Observe that \cite{npr} considers general convex domains. 
It may be instructive to read section 2.4 in \cite{wigmanr}. 
When dealing with the asymptotic distribution of the nodal length, 
the case of $Q_n$ is simpler.

\section{Comparison with the real-valued $2$-dimensional case} \label{s:discussion}
We begin with a common feature of both cases. 
As said above, in both cases, the second component in the Wiener expansion vanishes. 
This cancellation avoids having a higher order of the variances since the covariance functions are not in $L^2$.
\smallskip

Regarding the growing domain regime, that is, the setting of this paper. 
In view of the results in \cite{del}, 
theorem \ref{t:main} implies that the variance, and the expectation, 
of the nodal length are exactly of the same 
order as the volume of $B_R$. 
Recall that in the real-valued $2$-dimensional case the ratio variance/expectation 
was of the order of $\log({\rm area}(B))$ as the domain $B\uparrow\R^2$. 
\smallskip

This fact is connected to another important difference with the real-valued $2$-dimensional case. 
Here, each component of the Wiener expansion of $\cL(B)$ plays a role in the asymptotic variance, 
that is, there is no dominating term as in the real-valued $2$-dimensional case. 
This is due to the fact that the covariance function belongs to $L^4(\R^3)$ 
implying that all the integrals involved in the computation of the variance are convergent. 
To show this point, consider one of the terms involved in the variance of the sixth component of the expansion. 
\begin{equation*}
 \int_{B_R\times B_R}r(|x-y|)^6dxdy 
 	\sim \int_{B_R} \int_{B_R} r(|x-y|)^6 dxdy
	\sim \frac{2\pi^3 R^3}{3},
\end{equation*}
which is of the order of the volume of $B_R$. 
See section \ref{s:anciliary} for the computation of this sort of integrals. 
In the  real-valued $2$-dimensional case the analogous integral was still convergent but negligible w.r.t. to the variance of 
$I_4(B)$ since the covariance 
was in $L^6(R^2)\setminus L^4(\R^2)$. 
The existence of a single dominant chaotic component makes it possible to get quantitative central limit theorems using well-known 
bounds for the (e.g.: Kolmogorov, Wasserstein) distance between the distributions of r.v.s lying in a fixed chaotic component. We are 
not pursuing this kind of results in 
the present case where an infinite number of components contribute to the limit.
\medskip

Now, let us consider the high energy framework on a fixed domain. 
We follow remark 1.3. in \cite{npr}. 
Recall that for $k\geq 1$, we defined $\psi_k(x) = \psi(kx)$. 
Besides, we can relate the nodal lengths $\cL_k(B)$ of $\psi_k$ and $\cL(B)$ of $\psi$ by using their integral representation 
(see (2.23) in \cite{npr})
\begin{equation*}
 \cL_k(B) = \int_B \delta_0(\psi_k(x)) |\nabla\psi_k(x)|dx 
 = \frac1k\int_{k\cdot B} \delta_0(\psi(y)) |\nabla\psi(y)| dy 
 = \frac{\cL(k\cdot B)}{k}.
\end{equation*}
Here, as usual, we use $\delta_0$ as a shorthand notation for an approximation of the unity. To translate the former integral into the latter, we used the change of variables $y=kx$ and denoted $k\cdot B=\{kb:b\in B\}$. 
Hence, for fixed $B$, as $k\to\infty$, we have 
\begin{align*}
 \E[\cL_k(B)] &= \frac{\E[\cL(k\cdot B)]}{k} \sim \cs k^2;\\
 \Var[\cL_k(B)] &= \frac{\Var[\cL(k\cdot B)]}{k^2} \sim \cs k.
\end{align*}
We see, as above, that the ratio variance / square of the expectation tends to zero in the high energy limit and, thus, 
there is concentration of the probability 
in the sense that $\frac{\cL_k}{\E(\cL_k)}$ converges in probability to $1$ as $k\to\infty$. 
This ratio is $\approx 1/k^3$ in the complex-valued $3$-dimensional case and $\approx \log(k)/k^2$ in the real-valued $2$-dimensional one.

\section{Preliminaries} \label{s:prelim}
The departure point of the proof of theorem \eqref{t:main} is the Wiener expansion of $\cL(B_R)$ which we borrow from proposition 2 in  \cite{del}. 
With this expansion at hand, we bound the variance from below by that of one specific term (the fourth one) 
and compute it explicitly until we are sure of its positivity. 
\medskip

\noindent{\bf Hermite polynomials :} The building blocks of the Wiener expansion are Hermite polynomials which we define recursively by $H_0(x)=1$, 
$H_1(x)=x$ for $x\in\R$ and for $n\geq 2$ by
\[
H_n(x)=xH_{n-1}(x)-(n-1)H_{n-2}(x),\quad x\in\R.
\]
We need also the multi-dimensional Hermite polynomials :
\[
\tilde{H}_{\boldsymbol\alpha}({\boldsymbol y})=\prod^m_{i=1}H_{\alpha_i}(y_i),\quad 
{\boldsymbol \alpha}=(\alpha_i)_i\in\N^m\; \rm{and}\ {\boldsymbol y}=(y_i)_i\in\R^m.
\]
It is well known, see e.g. section 8.1. in \cite{pt}, that Hermite polynomials form a complete orthogonal system of $L^2(\phi_m(d{\boldsymbol 
y}))$ 
being $\phi_m$ the standard normal density function in $\R^m$. 
Hence, for $f\in L^2(\phi_m(d{\boldsymbol y}))$, 
we can write in the $L^2$-sense 
\begin{equation*}
f({\boldsymbol y})
=\sum^\infty_{q=0}\sum_{\boldsymbol \alpha \in \N^m,\,|{\boldsymbol \alpha}|=q}
f_{{\boldsymbol \alpha}}\tilde{H}_{\boldsymbol \alpha}({\boldsymbol y}),\quad {\boldsymbol y}\in \R^m,
\end{equation*}
with $|\balpha|=\sum^m_{i=1}\alpha_i$ and 
\begin{equation*}
f_{\boldsymbol \alpha}=\frac{1}{\balpha !}
\int_{\R^m}f({\boldsymbol y})\tilde{H}_{\boldsymbol \alpha}({\boldsymbol y})\phi_m(d{\boldsymbol y}),
\end{equation*}
with $\balpha !=\prod^m_{i=1}\alpha_i!$. 
We refer to $f_{\boldsymbol\alpha}$ as the $\boldsymbol\alpha$-th Hermite coefficient of $f$.
\medskip

\noindent{\bf Wiener expansion of $\cL(B_R)$ :} Let us introduce the Wiener expansion of the nodal length. 
Denote $\xi'(x)=(\xi_1(x),\xi_2(x),\xi_3(x))$ for the gradient of  $\xi(x)$ and $\eta'(x)=(\eta_1(x),\eta_2(x),\eta_3(x))$ for that of $\eta(x)$.  
Set
\begin{equation*}
Y(x) := 
\big(\xi(x),\eta(x),\sqrt{3}\xi'(x),\sqrt{3}\eta'(x)\big),
\end{equation*}
which is a standard normal random vector in $\R^8$ for each $x\in\R^3$. 
Finally, we need the formal Hermite coefficients of the Delta function \cite[Eq.5]{kl97}
\begin{equation}\label{e:b-alpha}
b_\alpha=\frac{1}{\alpha!\sqrt{2\pi}}H_\alpha(0), \quad \alpha\in\N. 
\end{equation}

Therefore, proposition 7 in \cite{del} states that  
\begin{align} \label{e:exp}
\tilde{\cL} (B_R) 
:= \frac{\cL(B_R)-\E[\cL(B_R)]} {\vol(B_R)} 
&=\frac{1}{3} \sum_{q\ge 1} I_{2q}(B_R), \\
I_{2q}(B_R) &= \sum_{{\boldsymbol \alpha}\in\N^8,\; |{\boldsymbol \alpha}|=2q} 
c_{\boldsymbol \alpha} \int_{B_R} \tilde{H}_{\boldsymbol \alpha}(Y(x)) dx, \notag
\end{align}
in the $L^2$-sense. 
The coefficients 
$c_{\boldsymbol \alpha}$ are defined by 
\[
c_{\boldsymbol \alpha} = b_{\alpha_1}b_{\alpha_2}a_{(\alpha_3,\dots,\alpha_8)},
\]
with $a_{(\alpha_3,\dots,\alpha_8)}$ the Hermite coefficient of $f({\boldsymbol y}) = \det{^{\bot}}({\boldsymbol y}) :=\det({\boldsymbol 
y}{\boldsymbol y}^\top)^{1/2}$, ${\boldsymbol y}\in\R^6$.
\medskip

\noindent{\bf Expectations of products of Hermite polynomials :} 
It is well known that for jointly Gaussian $X,Y$ 
\begin{equation} \label{e:mehler}
 \E[H_p(X)H_q(XY)] = \delta_{pq} p! \E[XY]^p, 
\end{equation}
where $\delta_{pq}$ is Kronecker's delta. 

For future use we recall a few more useful formulas, see e.g. \cite[Eq.6.76]{npr}. 
Let $X_1,X_2,X_3,X_4$ be standard normal r.v. with $\E(X_1X_2)=\E(X_3X_4)=0$, 
then
\begin{align} \label{e:diagrama}
 \E[H_2(X_1)H_2(X_2)H_2(X_3)H_2(X_4)] 
  &= 4\E[X_1X_3]^2\E[X_2X_4]^2 + 4\E[X_1X_4]^2\E[X_2X_3]^2 \notag\\
  &\qquad + 16\E[X_1X_3]\E[X_1X_4]\E[X_2X_3]\E[X_2X_4]; \notag \\
 \E[H_2(X_1)H_2(X_2)H_4(X_3)] 
  &= 24\E[X_1X_3]^2\E[X_2X_3]^2;\\ 
 \E[X_1X_2X_3X_4] &= \E[X_1X_3]\E[X_2X_4] + \E[X_1X_4]\E[X_2X_3]. \notag
\end{align}

\section{Proof of theorem \ref{t:main}} \label{s:proof}
Now, we can proceed to the proof of the main result.
\smallskip

\noindent{\bf Lower bound for the variance :} 
Since the components $I_{2q}(B_R):q\geq 1$ are orthogonal, 
we get
\begin{equation*} 
\Var(\cL(B_R)) = \frac13 \sum^\infty_{q=1}\Var(I_{2q}(B_R)) 
\geq \frac13 \Var(I_{4}(B_R)).
\end{equation*}
Thus, we reduce our attention to the fourth component $I_4(B_R)$. 
\medskip

\noindent{\bf Fourth component :} 
We can give the explicit expression of $I_4(B_R)$. 
From \eqref{e:exp} we get
\begin{align*}
I_4(B_R) &= b_0b_4a_0\int_{B_R}H_4(\xi(x))dx 
	   + b_0b_4a_0\int_{B_R}H_4(\eta(x))dx
		+ b^2_2a_0 \int_{B_R} H_2(\xi(x))H_2(\eta(x))dx\\ 
	&\qquad	+ \sum^3_{k=1} b^2_0a_{4e_k}\int_{B_R} H_4(\bar\xi_k(x))dx 
		+\sum^3_{k=1}b^2_0a_{4e_k}\int_{B_R} H_4(\bar\eta_k(x))dx  \\ 
	&\qquad+ \sum^3_{k=1}b_0b_2a_{2e_k}\int_{B_R} H_2(\xi(x))H_2(\bar\xi_k(x))dx  
		+ \sum^3_{k=1}b_0b_2a_{2e_k}\int_{B_R} H_2(\eta(x))H_2(\bar\eta_k(x))dx  \\
		&\qquad+ \sum^3_{k=1}b_0b_2a_{2e_k}\int_{B_R} H_2(\xi(x))H_2(\bar\eta_k(x))dx  
		+ \sum^3_{k=1}b_0b_2a_{2e_k}\int_{B_R} H_2(\eta(x))H_2(\bar\xi_k(x))dx  \\
	&\qquad+ \sum_{i\neq j}b^2_0a_{2e_i+2e_j}\int_{B_R} H_2(\bar\xi_i(x))H_2(\bar\xi_j(x))dx 
		  +\sum_{i\neq j}b^2_0a_{2e_i+2e_j}\int_{B_R} H_2(\bar\eta_i(x))H_2(\bar\eta_j(x))dx\\
	&\qquad	  +\sum_{i\neq j}b^2_0a_{2e_i+2e_j}\int_{B_R} H_2(\bar\xi_i(x))H_2(\bar\eta_j(x))dx.
\end{align*}
Here $\bar\xi_k=\sqrt{3}\xi_k$ and $\bar\eta_k=\sqrt{3}\eta_k$ 
and $\{e_1,e_2,\dots,e_6\}$ is the canonical basis in $\R^6$. 
Equation \eqref{e:b-alpha} gives
\[
b_0=\frac{1}{\sqrt{2\pi}}; b_2=-\frac{1}{2\sqrt{2\pi}}; 
b_4=\frac{1}{8\sqrt{2\pi}}. 
\]
The next lemma, 
whose proof is presented in section \ref{s:anciliary}, 
provides the explicit values of the relevant Hermite coefficients of $\det^\bot(\cdot)$.
\begin{lemma} \label{l:coef-a}
The first Hermite coefficients of $\det^\bot(\cdot)$ are :
\[
a_0 = 1;\ a_{2e_k} = \frac13;\ a_{2e_i+2e_j} = \frac19\ (i\neq j);\  
a_{4e_k} = -\frac59.
\]
\end{lemma}

Using the explicit values of the coefficients we have
\begin{equation*}
2\pi\cdot I_4 (B_R) =  A_1 + A_2 + A_3, 
\end{equation*}
with
\begin{align*}
A_1 &:= \frac{1}{8}\int_{B_R}H_4(\xi) + \frac{1}{8}\int_{B_R}H_4(\eta) 
	 +\frac{1}{4} \int_{B_R} H_2(\xi)H_2(\eta)
	 -\frac{5}{9}\sum^3_{k=1}\int_{B_R}H_4(\bar\xi_k) 
		- \frac{5}{9}\sum^3_{k=1}\int_{B_R}H_4(\bar\eta_k);  \\
A_2	&:=- \frac{1}{6} \sum^3_{k=1}\int_{B_R}H_2(\xi)H_2(\bar\xi_k)  
		- \frac{1}{6} \sum^3_{k=1}\int_{B_R}H_2(\eta)H_2(\bar\eta_k)  \\
		&\qquad- \frac{1}{6} \sum^3_{k=1}\int_{B_R}H_2(\xi)H_2(\bar\eta_k)  
		- \frac{1}{6}\sum^3_{k=1}\int_{B_R}H_2(\eta)H_2(\bar\xi_k);  \\
A_3	&:= \sum_{i\neq j}\frac{1}{9}\int_{B_R}H_2(\bar\xi_i)H_2(\bar\xi_j) 
		  +\sum_{i\neq j}\frac{1}{9}\int_{B_R}H_2(\bar\eta_i)H_2(\bar\eta_j)
		  +\sum_{i\neq j}\frac{1}{9}\int_{B_R}H_2(\bar\xi_i)H_2(\bar\eta_j).
\end{align*}
It follows that
\begin{equation*}
 \Var\big(2\pi\cdot I_4 (B_R)\big) =\sum^3_{i,j=1} \Cov(A_i,A_j).
\end{equation*}
This is a long computation.  
This partition, which may seem strange at first sight, 
will allow us to avoid considering quite a lot of terms.
\bigskip

Now, we state the equivalences for the variances and covariances of these terms, except for the variance of $A_2$ and $A_3$ which are certainly non negative and 
turn out to be innecssary.
\smallskip

Denote $A_1=\sum^5_{i=1}A_{1i}$. 
Then, as $\xi$ and $\eta$ are independent processes we have: 
\begin{align*}
 \Var(A_1) &= 2\Var(A_{11}) + \Var(A_{13}) +2\Var(A_{14})  
 + 4\Cov(A_{11},A_{14})\\
 &\sim 
  \frac{4\pi^3 R^3}{3} \frac{4362}{35}.
 \end{align*}

Similarly, denote $A_2=\sum^4_{i=1}A_{2i}$ and $A_3=\sum^3_{i=1}A_{3i}$. 
Then, 
\begin{align*}
\Cov(A_1,A_2) &= 
2\Cov(A_{11},A_{21}) + 2\Cov(A_{13},A_{23}) 
+ 2\Cov(A_{14},A_{21}) \\
& \gtrsim - \frac{4\pi^3 R^3}{3} \frac{4}{3}. \\
\Cov(A_1,A_3) &= 
2\Cov(A_{11},A_{31}) +\Cov(A_{13},A_{33}) + 2\Cov(A_{14},A_{31}) \\
& \gtrsim 
-\frac{4\pi^3 R^3}{3} \frac{1184}{105}. 
\end{align*}
and 
\begin{align*}
\Cov(A_2,A_3) &= 
2\Cov(A_{21},A_{31}) + 2\Cov(A_{23},A_{33}) \\
& \gtrsim 
-\frac{4\pi^3 R^3}{3} \frac{824}{525}.
\end{align*}
\smallskip

Therefore, since $\Var(A_2) \geq 0$ and $\Var(A_3) \geq 0$, we have
\begin{equation*}
\Var\big(I_4(B_R)\big) \gtrsim \frac{4\pi R^3}{3} \frac{7691}{350}. 
\end{equation*}
The result follows. 
The proofs of these relations are presented in section \ref{s:anciliary}.

\section{Ancillary computations} \label{s:anciliary}
In this section we include part of the long and 
sometimes tedious computations involved in the asymptotic variance of $I_4(B_R)$. 
The omitted computations are analogous to those which are presented here. 

\subsection{The coefficients}
We prove Lemma \ref{l:coef-a}. 
Note that $\det^\bot({\boldsymbol y}) = |{\boldsymbol y}\wedge{\boldsymbol y}|$. Here, $\wedge$ is the standard wedge product in $\R^3$.  
Assume that $Z_1 = (Z_{11},Z_{12},Z_{13})$ and $Z_2 = (Z_{21},Z_{22},Z_{23})$ are independent standard Gaussian r.v. in $\R^3$.
\begin{equation*}
 a_0 = \E (|Z_1\wedge Z_2|) = \E(\vol\{aZ_1+bZ_2:0\le a,b\le 1\})
	= \sqrt{2}\frac{\Gamma(1)}{\Gamma(1/2)}\cdot \sqrt{2}\frac{\Gamma(3/2)}{\Gamma(1)} 	=1.
\end{equation*}
Here we used the notation and results of \cite[p.305]{aw}  
for the expected volume of a random paralellepiped.
\smallskip

Since the distribution of $(Z_1,Z_2)$ is invariant under permutations of its coordinates, 
we deduce that $a_{2e_k}:k=1,\dots,6$ coincide. 
Besides, as $H_2(x)=x^2-1$,  
\begin{equation*}
 a_{2e_1} = \E (|Z_1\wedge Z_2| H_2(Z_{11})) 
 	= \frac13 \E (|Z_1\wedge Z_2| |Z_{1}|^2) - a_0.
\end{equation*} 
To connect the latter expectation with $a_0$ we use spherical coordinates, writing $z_1=\rho u$ with $\rho>0$ and $u\in S^2$. 
\begin{equation*} 
 \E (|Z_1\wedge Z_2| |Z_{1}|^2) = \int_{\R^3\times\R^3} |z_1\wedge z_2| |z_1|^2 
 	\phi_3(dz_1)\phi_3(dz_2) 
  	= \int_{\R^3} \phi_3(dz_2) \int_{S^2} |u\wedge z_2| du \int^\infty_{0} 
	  	\rho^2\rho^2\rho\phi(d\rho).
\end{equation*} 
Here we have a $\rho^2$ from the Jacobian, another $\rho^2$ from the squared norm and the final $\rho$ from 
the wedge product. 
In the same way, we get
\begin{equation*}
 a_{0} = \int_{\R^3\times\R^3} |z_1\wedge z_2| \phi_3(dz_1)\phi_3(dz_2)
      = \int_{\R^3} \phi_3(dz_2) \int_{S^2} |u\wedge z_2| du \int^\infty_{0} 	  	\rho^2\rho\phi(d\rho).
\end{equation*} 
Hence, denoting $m_k:=\int^\infty_{0} \rho^{k}\phi(d\rho)$ and using the recurrence 
$m_{k+1}=km_{k-1}$, we have
\begin{equation} \label{e:A}
 \E (|Z_1\wedge Z_2| |Z_{1}|^2) = \frac{m_5}{m_3}a_0 = 4. 
\end{equation} 
Thus,
\begin{align*}
 a_{2e_1} = \frac13 \frac{m_5}{m_3}a_0 -a_0 = \frac13 4-1 = \frac13. 
\end{align*}
\smallskip
 
Again, for $i\neq j$, all the $a_{2e_i+2e_j}$ coincide. 
Besides,
\begin{equation*}
 a_{2e_1+2e_4} = \E (|Z_1\wedge Z_2| H_2(Z_{11})H_2(Z_{21})) 
 	=  \E (|Z_1\wedge Z_2| Z^2_{11} Z^2_{21} ) - a_{2e_1} - a_{2e_4} -a_0.
\end{equation*} 
Here we used that 
$H_2(x)H_2(y) = x^2y^2 - H_2(x) - H_2(y) - 1$. 
Thus 
\begin{multline} \label{e:B}
  \E (|Z_1\wedge Z_2| Z^2_{11} Z^2_{21} ) 
  	= \frac19 \E (|Z_1\wedge Z_2| |Z_{1}|^2 |Z_{2}|^2 ) \\ 
	= \frac19 \int_{S^2\times S^2} |u_1\wedge u_2|du_1du_2 
		\int^\infty_0 \rho^2_1\rho^2_1\rho_1\phi(d\rho_1)
		\int^\infty_0 \rho^2_2\rho^2_2\rho_2\phi(d\rho_2)  
	= \frac19 \frac{m^2_5}{m^2_3} a_0 = \frac{16}{9},
\end{multline}
where we used spherical coordinates in both $z_1$ and $z_2$
and the fact that permuting the indexes does not changes the expectation. 
Thus
\begin{equation*}
 a_{2e_1+2e_4} = \frac19 \frac{m^2_5}{m^2_3} a_0 - \frac23 - 1 
 	= \frac19.
\end{equation*}
\smallskip

Finally, $a_{4e_k}:k=1,\dots,6$ coincide. 
Besides,
\begin{align*}
  a_{4e_k} &= \E (|Z_1\wedge Z_2| H_4(Z_1)) 
  		= \E ( |Z_1\wedge Z_2| (Z^4_1 -6Z^2_1 +3) )\\
		&= \frac13 \E (|Z_1\wedge Z_2| |Z_1|^4) -2 \E (|Z_1\wedge Z_2| Z^2_{11} Z^2_{21} ) 
		-2 \E (|Z_1\wedge Z_2| |Z_1|^2) + 3 a_0
\end{align*}
Here we used that $|(a_1,a_2,a_3)|^4 = a^4_1+a^4_2+a^4_3 
+ 2(a^2_1a^2_2+a^2_1a^2_3+a^2_2a^2_3)$.

The second and third terms have been computed in \eqref{e:A} and \eqref{e:B}. 
We look at the first one.
\begin{align*}
  \E (|Z_1\wedge Z_2| |Z_1|^4) &= \int_{\R^3\times \R^3}|z_1\wedge z_2| |z_1|^4   \phi_3(dz_1) \phi_3(dz_2) \\
  	&= \int_{\R^3}  \phi_3(dz_2) \int_{S^2} |u\wedge z_2|du \int^\infty_0 \rho^2 \rho^4 \rho 
	\phi(d\rho) 
	= \frac{m_7}{m_3}a_0
	= 24.
\end{align*}
Therefore,
\begin{equation*}
  a_{4e_k} = \frac13 \frac{m_7}{m_3}a_0 - 2\frac19 \frac{m^2_5}{m^2_3}a_0 
  - 2\frac{m_5}{m_3}a_0+3a_0 
  	= -\frac59.
\end{equation*}

\subsection{The variances and covariances of $A_{ij}$} 
We use \cite{GR} and \cite{wolfram} to compute the radial integrals 
(i.e: those w.r.t. $\rho$).
\smallskip

Recall that $B_R$ is the ball of radius $R$ centered at $0$ 
and set $B_R(z) = B_R + z$.
Also recall that 
\begin{align*} 
 r(u) = \sinc(u) &= \frac{\sin(u)}{u}, \notag\\
 \sinc'(u) &= \frac{\cos(u)}{u} - \frac{\sinc(u)}{u} 
 	= \frac{u\cos(u)-\sin(u)}{u^2}; \notag\\
 \sinc''(u) &= -\sinc(u) + \frac{2\sinc(u)}{u^2} - \frac{2\cos(u)}{u^2}.
\end{align*}

Besides, denoting $r_{k0}(x,y)=\frac{\partial r(x,y)}{\partial x_k}$ 
(and using similar notations for $y$), we have 
\begin{equation} \label{e:rk}
r_{k0}(x,y) = \sinc'(|x-y|)\Delta_k, \qquad 
r_{0k}(x,y) = -\sinc'(|x-y|)\Delta_k,
\end{equation}
with
\begin{equation*} 
\Delta_k  = \frac{x_k-y_k}{|x-y|}.  
\end{equation*}
Besides,
\begin{align} \label{e:rkk}
 r_{kk}(x,y) &= \Big[ \frac{\sinc'(|x-y|)}{|x-y|} -\sinc''(|x-y|) \Big]
   \frac{(x_k-y_k)^2}{|x-y|^2}   -\frac{\sinc'(|x-y|)}{|x-y|} \notag\\
   &=: A(|x-y|) \Delta^2_k - B(|x-y|),
\end{align}
with
\begin{align} \label{e:abdelta}
 A(u) &= \frac{\sinc'(u)}{u} -\sinc''(u) 
			= \frac{u^2\sin(u)-3\sin(u)+3u\cos(u)}{u^3} \notag\\
B(\rho) &= \frac{\sinc'(u)}{u}  
			= \frac{u\cos(u)-\sin(u)}{u^3}.
\end{align}
For $k\neq k'$: 
\begin{align} \label{e:rkkp}
 r_{kk'}(x,y) &= 
 \Big[ \sinc''(|x-y|) - \frac{\sinc'(|x-y|)}{|x-y|}\Big]
 \frac{(y_{k'}-x_{k'})(x_k-y_k)}{|x-y|^2} \notag\\
   &= A(|x-y|) \Delta_k \Delta_{k'}
\end{align}

To deal with the integrals of the angular parts $\Delta_k$ 
we consider spherical coordinates in $\partial B_\rho(x)$: 
\begin{equation*}
\begin{cases}
 x_1-y_1 &= \rho\sin(\theta)\cos(\varphi),\\
 x_2-y_2 &= \rho\sin(\theta)\sin(\varphi),\\
 x_3-y_3 &= \rho\cos(\theta).
\end{cases}
\end{equation*}

In the sequel we use sistematically \eqref{e:mehler} and \eqref{e:diagrama}.
\medskip

\begin{enumerate}
 \item We consider $\Var(A_1)$ in detail. 
We have.
\begin{itemize}
\item $\Var(A_{11}) = \Var(A_{12}) 
\sim \frac{4\pi^3 R^3}{3} \frac{3}{4}$,  
$\Var(A_{13}) \sim \frac{4\pi^3 R^3}{3} \frac{1}{2}$.
Indeed, using \eqref{e:diagrama} we get
\begin{align*}
 \Var(A_{11}) &= \frac{1}{64}\int_{B_R\times B_R}
 	\E( H_4(\xi(x)) H_4(\xi(y))) dxdy
	= \frac{4!}{64} \int_{B_R\times B_R} r(|x-y|)^4 dxdy.\\
 \Var(A_{13}) &=  \frac{1}{16}\int_{B_R\times B_R}
 	\E( H_2(\xi(x)) H_2(\xi(y))) \E( H_2(\eta(x)) H_2(\eta(y)))dxdy\\
	&= \frac{4}{16} \int_{B_R\times B_R} r(|x-y|)^4 dxdy.
\end{align*}

To compute this kind of integral we begin using the area formula 
as in the first display in the proof of proposition 5.1 in \cite{npr}. 
We have
\begin{align*}
\int_{B_R\times B_R} r(|x-y|)^4dxdy
	= \int^{2R}_{0} d\rho \int_{B_R} dx 
	\int_{\partial B_\rho(x)\cap B_R} r(|x-y|)^4dy,
\end{align*}
where $\partial B_\rho(x)$ is the sphere of radius $\rho$ centered at $x$. 
Following \cite[eq.(5.64)]{npr} we split this integral.
\begin{align*}
\int_{B_R\times B_R} r(|x-y|)^4dxdy 	
	&=  \int^{R}_{0} d\rho \int_{B_{R-\rho}} dx 
	\int_{\partial B_{\rho}(x)} r(|x-y|)^4dy \\
	&\qquad+  \int^{R}_{0} d\rho \int_{B_R-B_{R-\rho}} dx 
	\int_{\partial B_{\rho}(x)\cap B_R} r(|x-y|)^4dy\\
	&\qquad+  \int^{2R}_{R} d\rho \int_{B_R} dx 
	\int_{\partial B_{R}(x)\cap B_R} r(|x-y|)^4dy.
\end{align*}

Consider the first integral. 
Note that for $y\in \partial B_{\rho}(x)$ we have $|x-y| = \rho$. 
Hence,
\begin{align*}
&\int^{R}_{0} d\rho \int_{B_{R-\rho}} dx \int_{\partial B_{\rho}(x)} r(|x-y|)^4dy 
= \int^{R}_{0} \vol (B_{R-\rho}) {\rm area}(\partial B_{\rho}(x)) r(\rho)^4d\rho \\
&\qquad  
= \frac{16\pi^2}{3}\int^{R}_{0}  (R-\rho)^3 \rho^2 \frac{\sin(\rho)^4}{\rho^4}d\rho = \frac{16\pi^2}{3}R^3
	\int^{R}_{0}  \big(1-\frac{\rho}{R}\big)^3 \frac{\sin(\rho)^4}{\rho^2}d\rho \\
&\qquad \sim_{R\to\infty} \frac{16\pi^2}{3}R^3
	\int^{\infty}_{0}   \frac{\sin(\rho)^4}{\rho^2}d\rho 
= 2\frac{4\pi^3 R^3}{3}.
\end{align*}

Consider the second integral.
\begin{align*}
&\int^{R}_{0} d\rho \int_{B_R\setminus B_{R-\rho}} dx 
	\int_{\partial B_{\rho}(x)\cap B_R} r(|x-y|)^4dy 
	\leq \frac{16\pi^2}{3} \int^R_0 \big(R^3-(R-\rho)^3\big) 
		\rho^2 r(\rho)^4d\rho \\
	&\qquad \leq \frac{16\pi^2}{3} \int^R_0 \big(3R - 
		\frac{3R^2}{\rho}+\rho\big) \sin(\rho)^4d\rho
	= O(R^2\log(R)) = o(R^3).
\end{align*}

Similarly, for the third integral we have
\begin{align*}
&\int^{2R}_{R} d\rho \int_{B_R} dx 
	\int_{\partial B_{R}(x)\cap B_R} r(|x-y|)^4dy 
	\leq \int^{2R}_{R} \vol (B_{R}) {\rm area}(\partial B_{\rho}(x)) r(\rho)^4d\rho  \\ 		
	 &\qquad \leq   \frac{16\pi^2R^3}{3} \int^{2R}_{R}  \frac{\sin(\rho)^4}{\rho^2}d\rho 
	\leq   \frac{16\pi^2R^3}{3} \int^{2R}_{R}  \frac{1}{\rho^2}d\rho 
	= \frac{8\pi^2R^2}{3} = o(R^3). 
\end{align*}
\smallskip

In conclusion
\begin{equation*}
\int_{B_R\times B_R} r(|x-y|)^4dxdy 
= \frac{8\pi^3R^3}{3} + O(R^2\log(R)) 
\sim 2 \frac{4\pi^3R^3}{3},
\end{equation*}
as claimed.
\medskip

\item $\Cov(A_{11},A_{14}) = \Cov(A_{12},A_{15})
\sim - \frac{4\pi^3 R^3}{3}\frac{21}{5}$. 
Indeed,
\begin{equation*}
\Cov(A_{11},A_{14}) = -\frac{5}{72} \sum^3_{k=1}
\int_{B_R\times B_R}\E( H_4(\xi(x)) H_4(\bar{\xi}_k(y))) dxdy. 
\end{equation*}
Recall that $\bar{\xi}_k(x) = \sqrt{3}\xi_k(x)$. 
\begin{align*}
 \int_{B_R\times B_R}\E( H_4(\xi(x)) H_4(\bar{\xi}_k(y))) dxdy  
	&= 9\cdot 4! \int_{B_R\times B_R} 
	r_{0k}(x,y)^4dxdy\\
	&= 9\cdot 4! \int^{2R}_{0} d\rho \int_{B_R} dx 
		\int_{\partial B_\rho(x)\cap B_R} r_{0k}(x,y)^4dy.
\end{align*}
Recall that $r_{0k}(x,y)$ was computed in \eqref{e:rk}. 
We split the integral as before, 
the dominant term is 
\begin{align*}
 &\sum^3_{k=1} \int^{R}_{0} d\rho \int_{B_{R-\rho}} dx 
		\int_{\partial B_\rho(x)} r_{0k}(x,y)^4dy\\
 &= \sum^3_{k=1} \int^{R}_{0} \vol(B_{R-\rho}) \sinc'(\rho)^4
 \rho^2d\rho
 \int_{[0,\pi]\times[0,2\pi]} \Delta^4_k \sin(\theta)
 d\theta d\varphi
\end{align*}
The radial part is similar to the previous one.
\begin{align*}
 \int^{R}_{0} \vol(B_{R-\rho}) \sinc'(\rho)^4 \rho^2 d\rho 
 &= \frac{4\pi R^3}{3} \int^{R}_{0}  \big[ 1- \frac{\rho}{R}\big]^3
 \sinc'(\rho)^4 \rho^2 d\rho \\
 &\sim
 \frac{4\pi R^3}{3} \int^{\infty}_{0}  \sinc'(\rho)^4 \rho^2 d\rho 
 = \frac{4\pi R^3}{3} \frac{7\pi}{60}.
\end{align*}
The angular part 
equals $\frac{12\pi}{5}$. 
Indeed,
\begin{align*}
&\sum^3_{k=1} \int_{[0,\pi]\times[0,2\pi]} 
\Delta^4_k \sin(\theta)d\theta d\varphi \\
&\qquad = \int_{[0,\pi]\times[0,2\pi]}\Big[ \sin^5(\theta)\cos^4(\varphi) 
 + \sin^5(\theta)\sin^4(\varphi)
 + \cos^4(\theta) \sin(\theta) \Big] d\theta d\varphi \\
 &\qquad = \frac{16}{15}\cdot \frac{3\pi}{4} 
 + \frac{16}{15}\cdot \frac{3\pi}{4} + \frac{2}{5}2\pi 
 = \frac{12\pi}{5}.
\end{align*}
In the sequel we omit the explicit computations of the angular parts since they are analogous.
(The remainder term behaves as in the previous case). 
In conclusion
\begin{equation*}
\sum^3_{k=1}\int_{B_R\times B_R}
	\E( H_4(\xi(x)) H_4(\bar{\xi}_k(y))) dxdy  
 	\sim 9\cdot 4! \frac{4\pi R^3}{3} \frac{7\pi}{60} \frac{12\pi}{5}
	= 	\frac{4\pi R^3}{3} \frac{1512}{25}.
\end{equation*}
Multiplying by $-\frac{5}{72}$ the result follows.
\medskip

\item $\Var(A_{14}) = \Var(A_{15})
=  \frac{488}{7}$.  
Indeed,
\begin{align*}
\Var(A_{14}) = \frac{25}{81} \sum^3_{k,k'=1} 
\int_{B_R\times B_R}\E( H_4(\bar{\xi}_k(x)) H_4(\bar{\xi}_k(y))) dxdy.
\end{align*}

We begin with the diagonal terms.
\begin{equation*}
 \int_{B_R\times B_R}\E( H_4(\bar{\xi}_k(x)) H_4(\bar{\xi}_k(y))) dxdy  
= 81\cdot 4!  \int_{B_R\times B_R} r_{kk}(x,y)^4dxdy.
\end{equation*}
Recall that $r_{kk}$ is given by \eqref{e:rkk} and \eqref{e:abdelta}. 
Thus,
\begin{align*}
 \int_{B_R\times B_R} r_{kk}(x,y)^4dxdy 
 	&\sim \frac{4\pi R^3}{3}\int^{R}_{0} \big[1-\frac{\rho}{R}\big]^3 
	\int_{\partial B_{\rho}(x)} r_{kk}(x,y)^4 dy d\rho\\
	&\sim \frac{4\pi R^3}{3}\int^{\infty}_{0}  
	\int_{\partial B_{\rho}(x)} r_{kk}(x,y)^4 dy d\rho.
\end{align*}
The fourth power gives rise to several terms that we treat separately. 
For the first one, using spherical coordinates as before, we have 
\begin{align*}
&\sum^3_{k=1}\int^{\infty}_{0} 
\int_{\partial B_{\rho}(x)} A(|x-y|)^4 \Delta^8_k dxdy 
	= \int^{\infty}_{0} A(\rho)^4\rho^2 d\rho 
	\cdot \int_{[0,\pi]\times[0,2\pi]}\Delta^8_k 
	\sin(\theta) d\theta d\varphi \\
	&\qquad = \frac{11\pi}{140}\frac{12\pi}{9} 
	=\frac{11\pi^2}{105}.
\end{align*}
\smallskip

Moving to the next terms, 
with similar arguments we have
\begin{align*} 
\sum^3_{k=1} 
\int^{\infty}_{0} \int_{\partial B_{\rho}(x)} A(\rho)^3 
	\Delta_k^6B(\rho) dy d\rho &= \frac{\pi}{70} \frac{12\pi}{7} 
	= \frac{6\pi^2}{245}\\
\sum^3_{k=1} \int^{\infty}_{0} 
\int_{\partial B_{\rho}(x)} A(\|x-y\|)^2 \Delta_k^4 B(\|x-y\|)^2 dx dy 
	&= \frac{2\pi}{315}\frac{12\pi}{5} = \frac{8\pi^2}{525},\\
\sum^3_{k=1} \int^{\infty}_{0} \int_{\partial B_{\rho}(x)} A(\rho)
      \Delta_k^2 B(\rho)^3 dy d\rho
      &= \frac{17\pi}{3780}\frac{12\pi}{3} =\frac{17\pi^2}{945},\\
\int^{\infty}_{0} \int_{\partial B_{\rho}(x)} B(\rho)^4 dy d\rho
      &= 4\pi \frac{17\pi}{2835} =\frac{68\pi^2}{2835}
\end{align*}

Gathering all together.
\begin{align*}
 &\sum^3_{k=1}
 \int_{B_R\times B_R}\E( H_4(\bar{\xi}_k(x)) H_4(\bar{\xi}_k(y))) dxdy \\
 &\qquad \sim 81\cdot 4! \frac{4\pi R^3}{3} 
 \Big[ \frac{11\pi^2}{105} 
 -4 \frac{6\pi^2}{245} 
 + 6 \frac{8\pi^2}{525}
 -4 \frac{17\pi^2}{945}
 + 3 \frac{68\pi^2}{2835}  \Big] \\
 &\qquad = 81\cdot 4! \frac{4\pi^3 R^3}{3} \frac{361}{3675}
 = \frac{4\pi^3 R^3}{3} \frac{233928}{1225}.
 \end{align*}
 \end{itemize}

We move to the case $k\neq k'$. 
\begin{equation*}
 \int_{B_R\times B_R}\E( H_4(\bar{\xi}_k(x)) H_4(\bar{\xi}_{k'}(y))) dxdy  
= 81\cdot 4! \int_{B_R\times B_R} r_{kk'}(x,y)^4dxdy.
\end{equation*}
Recall that $r_{kk'}$, $k\neq k'$, is given by \eqref{e:rkkp}. 
\begin{align*}
&\sum_{k\neq k'}\int_{B_R\times B_R} r_{kk'}(x,y)^4dxdy 
	\sim \frac{4\pi R^3}{3}\sum_{k\neq k'} \int^{\infty}_{0}  
	\int_{\partial B_{\rho}(x)} r_{kk'}(x,y)^4 dy d\rho\\
	&\qquad= \frac{4\pi R^3}{3}\int^{\infty}_{0} 
	A(\rho)^4\rho^2 d\rho\cdot \sum_{k\neq k'} \int_{[0,\pi]\times[0,2\pi]} 
\Delta^4_k \Delta^4_{k'} \sin(\theta)d\theta d\varphi   \\
	&\qquad = \frac{4\pi R^3}{3} \frac{11\pi}{140} 
	\frac{24\pi}{105}  = \frac{4\pi^3 R^3}{3} \frac{22}{1225}. 
\end{align*}
Thus, multiplying this last result by $81\cdot 4!$ and
summing it to the diagonal terms:
\begin{equation*}
\sum^3_{k, k'=1}\int_{B_R\times B_R}\E( H_4(\bar{\xi}_k(x)) H_4(\bar{\xi}_{k'}(y))) dxdy  
=  \frac{4\pi^3 R^3}{3} \Big[ \frac{233928}{1225} 
  + \frac{42768}{1225} \Big]
  = \frac{4\pi^3 R^3}{3} \frac{39528}{175}.
\end{equation*}
Multiplying by $25/81$ we get the result. 
\medskip

\item Let us consider $\Cov(A_{1},A_{2})$. 
There are many terms which vanish since $\xi$ and $\eta$ are independent random fields. 
\begin{multline*}
 \Cov(A_{11},A_{23}) = \Cov(A_{11},A_{24}) 
= \Cov(A_{12},A_{23}) = \Cov(A_{12},A_{24}) =\Cov(A_{13},A_{21}) \\
= \Cov(A_{13},A_{24}) =\Cov(A_{14},A_{22}) = \Cov(A_{14},A_{23}) 
= \Cov(A_{14},A_{24})\\ 
= \Cov(A_{15},A_{21}) = \Cov(A_{15},A_{23}) 
= \Cov(A_{15},A_{24}) = 0.
\end{multline*}

We treat the rest of them. 
\begin{itemize}
\item $\Cov(A_{11},A_{21}) = \Cov(A_{12},A_{22})
\sim  -\frac{4\pi^3 R^3}{3} \frac12$. 
Indeed,
\begin{multline*}
 \sum^3_{k=1}\int_{B_R\times B_R} 
 \E \big[ H_4(\xi(x)) H_2(\xi(y)) H_2(\bar{\xi}_{k}(y)) \big] dxdy \\
 = 3\cdot 4! \sum^3_{k=1}  \int_{B_R\times B_R} r(x,y)^2 r_{0k}(x,y)^2dxdy\\
 = 3\cdot 4! \sum^3_{k=1}  \int_{B_R\times B_R} 
 \sinc(\rho)^2 \sinc'(\rho)^2 \Delta^2_kdxdy
 \sim 3\cdot 4!\frac{4\pi R^3}{3} 
 \frac{\pi}{12} 4\pi  =  \frac{4\pi^3 R^3}{3}  24.
\end{multline*}
Multiplying by $-1/48$ we get the result.
\medskip

\item $\Cov(A_{13},A_{23}) = \Cov(A_{13},A_{24})
\sim  -\frac{4\pi^3 R^3}{3} \frac{1}{6}$. Indeed,
\begin{align*}
 &\sum^3_{k=1}\int_{B_R\times B_R} 
 \E \big[ H_2(\xi(x)) H_2(\xi(y)) \big] 
 \E \big[ H_2(\eta(y)) H_2(\bar{\eta}_k(y)) \big] dxdy\\
 &\qquad= 3\cdot 4 \sum^3_{k=1} \int_{B_R\times B_R} 
 r(x,y)^2 r_{0k}(x,y)^2dxdy
 \sim 3\cdot 4 \frac{4\pi R^3}{3} \frac{\pi^2}{3}
 =\frac{4\pi^3 R^3}{3} 4.
\end{align*}
We used the computations in $\Cov(A_{11},A_{21})$.
Multiplying by $-1/24$ we get the result.
\medskip

\item $\Cov(A_{14},A_{21}) = \Cov(A_{15},A_{22}) \geq 0$.
We have
\begin{multline*}
 \Cov(A_{14},A_{21}) = 
 \frac{5}{54} \sum^3_{k,k'=1}\int_{B_R\times B_R} 
 \E \big[ H_4(\bar{\xi}_k(x)) H_2(\xi(y))  H_2(\bar{\xi}_{k'}(y)) \big] dxdy\\
 = \frac{5}{54}  27\cdot 4! \sum_{k,k'=1} \int_{B_R\times B_R} 
 r_{k0}(x,y)^2 r_{kk'}(x,y)^2dxdy \geq 0.
\end{multline*}
Thus, we can omit it.
\end{itemize}
\medskip

\item We move to $\Cov(A_{1},A_{3})$. 
Again there are several vanishing terms. 
\begin{multline*}
\Cov(A_{11},A_{32}) = \Cov(A_{11},A_{33}) 
= \Cov(A_{12},A_{31}) = \Cov(A_{12},A_{33}) =\Cov(A_{13},A_{31}) \\
= \Cov(A_{13},A_{32}) = \Cov(A_{14},A_{32}) = \Cov(A_{14},A_{33}) 
= \Cov(A_{15},A_{31}) = \Cov(A_{15},A_{33}) = 0.
\end{multline*}

We look at the rest of the terms.
\begin{itemize}

\item $\Cov(A_{11},A_{31}) = \Cov(A_{12},A_{32}) \geq 0$. 
Indeed,
\begin{multline*}
 \Cov(A_{11},A_{31}) = \frac{1}{72} \sum_{i\neq j}\int_{B_R\times B_R} 
 \E \big[ H_4(\xi(x)) H_2(\bar{\xi}_i(y)) H_2(\bar{\xi}_{j}(y)) \big] dxdy \\
 = \frac{1}{72} 9\cdot 4!  \sum_{i\neq j} \int_{B_R\times B_R} 
 r_{0i}(x,y)^2 r_{0j}(x,y)^2dxdy \geq 0. 
\end{multline*}
\smallskip

\item $\Cov(A_{13},A_{33}) \geq 0$.
we have,
\begin{multline*}
 \Cov(A_{13},A_{33}) = \frac{1}{36} \sum_{i\neq j}\int_{B_R\times B_R} 
 \E \big[ H_2(\xi(x)) H_2(\bar{\xi}_i(y)) \big] 
 \E \big[ H_2(\eta(x)) H_2(\bar{\eta}_j(y)) \big] dxdy\\
  = \frac{1}{36} 9\cdot 4 \sum_{i\neq j} \int_{B_R\times B_R} 
  r_{0i}(x,y)^2 r_{0j}(x,y)^2dxdy \geq 0.
\end{multline*}
\smallskip

\item $\Cov(A_{14},A_{31}) = \Cov(A_{15},A_{32})
\sim - \frac{4\pi^3R^3}{3} \frac{592}{105}$. 
We have,
\begin{multline*}
 \sum^3_{k=1} \sum_{i\neq j} \int_{B_R\times B_R} 
 \E \big[ H_4(\bar{\xi}_k(x)) H_2(\bar{\xi}_i(y)) H_2(\bar{\xi}_{j}(y)) \big] 
 dxdy \\
 = 81\cdot 4! \sum^3_{k=1}  \sum_{i\neq j}
 \int_{B_R\times B_R} r_{ki}(x,y)^2 r_{kj}(x,y)^2dxdy
\end{multline*}
\smallskip

We begin with the case $k\neq i, k\neq j$, $i\neq j$, 
we write $k\neq i\neq j$ for short. 
\begin{align*}
\sum_{k\neq i\neq j} \int_{B_R\times B_R} r_{ki}(x,y)^2 r_{kj}(x,y)^2dxdy 
&= \sum_{k\neq i\neq j}
\int_{B_R\times B_R} A(|x-y|)^4 \Delta^4_k \Delta^2_i \Delta^2_j dxdy\\
& = \frac{4\pi R^3}{3} \frac{11\pi}{140}
\frac{24\pi}{315}
= \frac{4\pi^3 R^3}{3} \frac{22}{3675}.
\end{align*}
Now, we consider the case $k = i\neq j$, 
($k = j \neq i$ is equal). 
\begin{align*}
&\sum_{k\neq j} \int_{B_R\times B_R} r_{kk}(x,y)^2 r_{kj}(x,y)^2dxdy\\ 
&= \sum_{k\neq j} \int_{B_R\times B_R} 
(A(|x-y|)^2 \Delta^2_k -B(|x-y|))^2 
A(|x-y|)^2 \Delta^2_k \Delta^2_j dxdy.
\end{align*}
The first integral
\begin{align*}
 \sum_{k\neq j} \int_{B_R\times B_R} A(|x-y|)^4 
\Delta^6_k \Delta^2_j dxdy
 &\sim  \frac{4\pi R^3}{3} \frac{11\pi}{140}  \frac{8\pi}{21}
= \frac{4\pi^3 R^3}{3} \frac{22}{735}.
\end{align*}
The radial part is as in the previous case. 
The second integral.
\begin{align*}
\sum_{k\neq j} \int_{B_R\times B_R} A(|x-y|)^3 B(|x-y|) 
\Delta^4_k \Delta^2_j dxdy
 \sim  \frac{4\pi R^3}{3} \frac{\pi}{70} \frac{24\pi}{35}
 = \frac{4\pi^3 R^3}{3} \frac{12}{1225}. 
\end{align*} 
The third integral is
\begin{align*}
 \sum_{k\neq j} \int_{B_R\times B_R} A(|x-y|)^2 B(|x-y|)^2 
\Delta^2_k \Delta^2_j dxdy
 &\sim  \frac{4\pi R^3}{3} \frac{2\pi}{315} \frac{8\pi}{5} 
 = \frac{4\pi^3 R^3}{3} \frac{16}{1575}.
\end{align*} 
\smallskip

Thus,
\begin{align*}
 &\sum^3_{k=1} \sum_{i\neq j} \int_{B_R\times B_R} 
 \E \big[ H_4(\bar{\xi}_k(x)) H_2(\bar{\xi}_i(y)) H_2(\bar{\xi}_{j}(y)) \big] 
 dxdy\\
&\qquad \sim 81\cdot 4! \frac{4\pi^3 R^3}{3} 
\Big[ \frac{22}{3675} + 2\big( \frac{22}{735} 
-2 \frac{12}{1225} + \frac{16}{1575}\big) \Big] 
= \frac{4\pi R^3}{3} \frac{15984}{175}.
\end{align*}
The result follows multiplying by $-5/81$.
\end{itemize}
\bigskip

\item Finally, we consider $\Cov(A_{2},A_{3})$.
We have
\begin{multline*}
\Cov(A_{21},A_{32}) = \Cov(A_{21},A_{33}) 
= \Cov(A_{22},A_{31}) = \Cov(A_{22},A_{33}) \\
= \Cov(A_{23},A_{31}) = \Cov(A_{23},A_{32})  
= \Cov(A_{24},A_{31}) = \Cov(A_{24},A_{32})  = 0.
\end{multline*}

Besides, we have.
\begin{itemize}

\item $\Cov(A_{21},A_{31}) = \Cov(A_{22},A_{32}) 
= -\frac{4\pi^3R^3}{3}\frac{1304}{3675} $.
\begin{align*}
 &\sum^3_{k=1}\sum_{i\neq j} \int_{B_R\times B_R} 
 \E \big[ H_2(\xi(x)) H_2 (\bar{\xi}_k(x)) 
 H_2(\bar{\xi}_i(y)) H_2(\bar{\xi}_{j}(y)) \big] dxdy\\
 &\qquad= 27\cdot 4 
 \sum^3_{k=1}\sum_{i\neq j} \int_{B_R\times B_R} 
 r_{0i}(x,y)^2 r_{kj}(x,y)^2dxdy\\
  &\qquad+ 27\cdot  4 
  \sum^3_{k=1}\sum_{i\neq j} \int_{B_R\times B_R} 
 r_{0j}(x,y)^2 r_{ki}(x,y)^2dxdy\\
  &\qquad+ 27\cdot 16
  \sum^3_{k=1}\sum_{i\neq j} \int_{B_R\times B_R} 
 r_{0i}(x,y) r_{0j}(x,y) r_{ki}(x,y) r_{kj}(x,y) dxdy
\end{align*}

We begin with the first two integrals 
(their sum over $k, i\neq j$ coincide). 
Consider $k\neq i$, $k\neq j$ (and $i\neq j$), 
we write $i\neq k\neq j$ for short. 
\begin{multline*}
\sum_{i\neq k\neq j} 
\int_{B_R\times B_R}  r_{0j}(x,y)^2 r_{ki}(x,y)^2dxdy \\
 = \sum_{i\neq k\neq j}
 \int_{B_R\times B_R} \sinc'(|x-y|)^2 A(|x-y|)^2
 \Delta^2_i \Delta^2_k \Delta^2_j dxdy  
  \sim  \frac{4\pi R^3}{3} \frac{23\pi}{420} 
 \frac{8\pi}{35} 
 = \frac{4\pi^3 R^3}{3} \frac{46}{3675}.
\end{multline*}

Now, take $k=j\neq i$. 
\begin{align*}
\sum_{i\neq k} 
\int_{B_R\times B_R}  r_{0k}(x,y)^2 r_{ki}(x,y)^2dxdy
&= \sum_{i\neq k}
 \int_{B_R\times B_R} \sinc'(|x-y|)^2 A(|x-y|)^2
 \Delta^2_i \Delta^4_k dxdy  \\
 & \sim \frac{4\pi R^3}{3} \frac{23\pi}{420} \frac{24\pi}{35}
 = \frac{4\pi^3 R^3}{3} \frac{46}{1225}. 
\end{align*}

Finally, take $k=i\neq j$.
\begin{align*}
&\sum_{j\neq k} 
\int_{B_R\times B_R} r_{0j}(x,y)^2 r_{kk}(x,y)^2dxdy\\
 &\qquad = \sum_{j\neq k}
 \int_{B_R\times B_R} \sinc'(|x-y|)^2 
 ( A(|x-y|) \Delta^2_k -B(|x-y|))^2
 \Delta^2_j  dxdy  \\
 &\qquad = \frac{4\pi R^3}{3} \Big[ 
 \frac{23\pi}{420} \frac{24\pi}{35} - 2 \frac{\pi}{42} 
 \frac{8\pi}{5} + \frac{2\pi}{105} 12\pi \Big]
 = \frac{4\pi^3 R^3}{3} \frac{698}{3675}.
\end{align*}
\medskip

We move to the third integral. 
Take $k\neq i$, $k\neq j$ and $i\neq j$. 
Write $i\neq k\neq j$ for short.
\begin{multline*}
\sum_{i\neq k\neq j} \int_{B_R\times B_R} 
 r_{0i}(x,y) r_{0j}(x,y) r_{ki}(x,y) r_{kj}(x,y) dxdy\\
 = \sum_{i\neq k\neq j} \int_{B_R\times B_R} 
 \sinc'(|x-y|)^4 \Delta^2_i\Delta^2_j \Delta^2_k\  dxdy
 \sim \frac{4\pi R^3}{3} \frac{7\pi}{60} \frac{8\pi}{35}
 = \frac{4\pi^3 R^3}{3}\frac{2}{75}.
\end{multline*}

Finally, take $k=j\neq i$ (the case $k=i\neq j$ yields the same value). 
\begin{align*}
& \sum_{i\neq k} \int_{B_R\times B_R} 
 r_{0i}(x,y) r_{0k}(x,y) r_{ki}(x,y) r_{kk}(x,y) dxdy\\
 &\qquad = \sum_{i\neq k} \int_{B_R\times B_R} 
 \sinc'(|x-y|)^2 A(|x-y|) (A(|x-y|)\Delta^2_k-B(|x-y|)) 
 \Delta^2_i \Delta^2_k\  dxdy\\
 &\qquad \sim \frac{4\pi R^3}{3} \Big[
 \frac{23\pi}{420} \frac{24\pi}{35} -  \frac{\pi}{42} \frac{8\pi}{5}\Big] 
 = - \frac{4\pi^3 R^3}{3} \frac{2}{3675}.
\end{align*}
\smallskip

Gathering all together
\begin{align*}
 &\sum^3_{k=1}\sum_{i\neq j}\int_{B_R\times B_R} 
 \E \big[ H_2(\xi(x)) H_2 (\bar{\xi}_k(x)) 
 H_2(\bar{\xi}_i(y)) H_2(\bar{\xi}_{j}(y)) \big] dxdy\\
 &\qquad= 27\cdot 4  \frac{4\pi^3 R^3}{3} 
 \Big[ 2 \big( \frac{46}{3675} +\frac{46}{1225}
 + \frac{698}{3675}\big) 
 + 4 \big( \frac{2}{75} -2\frac{2}{3675} \big) \Big] 
  = 27 \cdot 4  \frac{4\pi^3 R^3}{3} \frac{652}{3675}. 
\end{align*} 
Multiplying by $-1/54$ we get the result.
\medskip

\item $\Cov(A_{23},A_{33}) = \Cov(A_{24},A_{33})
\sim - \frac{4\pi^3 R^3}{3} \frac{316}{735}$.
\begin{align*}
 &\sum^3_{k=1}\sum_{i\neq j}\int_{B_R\times B_R} 
 \E \big[ H_2(\xi(x)) H_2(\bar{\xi}_i(y)) \big] 
 \E \big[ H_2(\bar{\eta}_k(x)) H_2(\bar{\eta}_j(y)) \big] dxdy\\
 &\qquad= 27\cdot 4 
 \sum^3_{k=1}\sum_{i\neq j} \int_{B_R\times B_R} 
 r_{0i}(x,y)^2 r_{kj}(x,y)^2dxdy
= 27\cdot 4  \frac{4\pi^3 R^3}{3} \frac{158}{735}. 
\end{align*}
Multiplying by $-1/54$ we get the result. 
This integral coincides with the first two terms in 
$\Cov(A_{21},A_{31})$.

\end{itemize}
\end{enumerate}

 
\end{document}